\documentclass[11pt]{article}
\usepackage{amsmath, amssymb, amsthm, hyperref, geometry,enumitem, comment}
\usepackage{enumitem}
\usepackage{tikz-cd}
\usepackage{url}
\usepackage{mathrsfs}
\usepackage{amsfonts}
\let\svthefootnote\thefootnote
\newcommand\freefootnote[1]{%
  \let\thefootnote\relax%
  \footnotetext{#1}%
  \let\thefootnote\svthefootnote%
}

\newtheorem{theorem}{Theorem}[section]
\newtheorem{lemma}[theorem]{Lemma}

\newtheorem{proposition}[theorem]{Proposition}
\theoremstyle{definition}
\newtheorem{definition}[theorem]{Definition}

\newtheorem*{pf}{Proof}
\newtheorem*{pf2}{Proof of Theorem 1.1}
\theoremstyle{remark}
\newtheorem{remark}[theorem]{Remark}

\begin{document}
\title{\normalsize\MakeUppercase{\textbf{On the Existence of Good Minimal Models for Kähler Varieties with Projective Albanese Map}}}
\author{\small\MakeUppercase{Yu-Ting Huang}}
\date{}

\maketitle
\freefootnote{\textit{Department of Mathematics, University of Utah, 155 South 1400 East, JWB 321, Salt Lake City, UT 84112, USA. E-mail address: ythuang@math.utah.edu}}

\begin{abstract}
In this article, we establish the existence of a good minimal model for a compact Kähler klt pair $(X, B)$ when the Albanese map of $X$ is a projective morphism and the general fiber of $(X, B)$ has a good minimal model.
\end{abstract}

\section{Introduction}

In this paper, we establish the following result.

\begin{theorem}\label{1}
Let $a_X: X\to A$ be the Albanese morphism of a compact Kähler variety $X$ where $(X, B)$ is a klt pair. Assume $a_X$ is projective with connected fibers. Let $F$ be the general fiber of $a_X: X\to a_{X}(X)$. Suppose $(F, B_F)$ has a good minimal model, then $(X, B)$ has a good minimal model.
\end{theorem}

We say the pair $(X, B)$ is a good minimal model over $U$ if $(X, B)$ is a minimal model and $K_X+B$ is semiample over $U$. The existence of good minimal models has been a central topic in higher dimensional geometry, closely related to the Minimal Model Conjecture and the Abundance Conjecture. The existence of good minimal models has been established in various settings. In \cite{fujino2009maximal}, Fujino first established that a smooth projective variety of maximal Albanese dimension possesses a good minimal model. This result was later generalized to projective varieties of arbitrary Albanese dimension in \cite{lai2011varieties}. Specifically, Lai proved the following theorem. 

\begin{theorem}[{\cite[Theorem 4.5]{lai2011varieties}}]\label{laimain}
Let $X$ be a smooth projective variety with Albanese map $\alpha: X\to A$. If the general fiber $F$ of $\alpha$ has dimension no more than $3$ with $\kappa(F)\geq 0$, then $X$ has a good minimal model.
\end{theorem}

The argument in Lai’s proof also verifies that if $F$ has a good minimal model, then $X$ must also have a good minimal model.

In recent years, the Minimal Model Program has been further developed for projective morphisms between compact Kähler varieties, see \cite{DAS2024109615} and \cite{fujino2022minimal}. However, generalizing \cite{fujino2009maximal} or Theorem~\ref{laimain} to a Kähler variety $X$ presents several technical challenges. 

For example, in an algebraic setting, following the method in \cite{fujino2009maximal}, it is clear that a relative minimal model over an Abelian variety is actually an absolute minimal model. Indeed, let $(X,B)$ be a projective klt pair, and $a: X\to A$ a projective morphism such that $K_X+B$ is nef over $A$. Suppose $X$ contains a $(K_X+B)$-negative rational curve $C$. Since $A$ contains no rational curve, $C$ must be contracted by $a$, which contradicts the assumption. Thus, by the cone theorem, $K_X+B$ is nef. However, the cone theorem remains unknown in the Kähler setting.

The proof of Theorem~\ref{laimain} is divided into three cases -- $\kappa(X) = -\infty$, $\kappa(X) = 0$, and $\kappa(X)>0$. Using generic vanishing techniques, it is shown that $\kappa(X)\geq 0$. Moreover, considering the Iitaka fibration, one shows that it suffices to settle the case $\kappa(X) = 0$ (cf. \cite[Theorem 4.4]{lai2011varieties}). Unfortunately, in the Kähler setting, the Iitaka fibration may not be a projective morphism, so \cite[Theorem 4.4]{lai2011varieties} cannot be directly generalized.

Note that Theorem~\ref{1} in the maximal Albanese dimension was proved by Das and Hacon.

\begin{theorem}[{\cite[Theorem 1]{das2024existence}}]
    Let $(X,B)$ be a compact Kähler klt pair of maximal Albanese dimension. Then $(X, B)$ has a good minimal model.
\end{theorem}

The rough idea of the proof of the theorem is as follows. Assuming that $X$ has maximal Albanese dimension, by \cite{DAS2024109615} and \cite{fujino2022minimal}, we may run a relative MMP over $A$, the Albanese variety of $X$, so that we may assume $K_X+B$ is nef over $A$. Since the cone theorem is not known in this context, the proof is instead concluded based on the approach in \cite{cao2020rational}.

In this article, we generalize Theorem~\ref{laimain} to compact Kähler klt pair. We will use \cite{DAS2024109615} and \cite{fujino2022minimal} to run a relative minimal model program. Then we prove Theorem~\ref{main1} to replace the cone theorem. Finally, we adapt the methods of \cite{lai2011varieties} and \cite{das2024existence} to complete the proof of Theorem~\ref{1}.

\section{Preliminaries}

\subsection{Minimal Model Program for Kähler Varieties}

In this subsection, we first recall the following result concerning the Minimal Model Program in the Kähler setting. For the definition of singularities of analytic varieties, refer to \cite{DAS2024109615} and \cite{fujino2022minimal}.

\begin{theorem}[cf. {\cite[Theorem 1.3]{DAS2024109615}} and \cite{fujino2022minimal}]\label{fg}
Let $f: (X,B)\to Y$ be a proper surjective morphism of compact analytic varieties where $(X, B)$ is a Kähler klt pair. Then the relative canonical $\mathcal{O}_Y$-algebra
$$R(X/Y, K_X+B):=\bigoplus_{m\geq 0}f_*\mathcal{O}_X(m(K_X+B))$$ is finitely generated.
\end{theorem}

\begin{theorem}[cf. {\cite[Theorem 1.4]{DAS2024109615}} and \cite{fujino2022minimal}]\label{kahlermmp}
    Let $\pi:(X, B)\to U$ be a projective morphism of compact normal analytic varieties where $X$ is a $\mathbb Q$-factorial (cf. \cite[Definition 2.7]{DAS2024109615}) and $(X, B)$ is a klt pair. Then,
    \begin{enumerate}[label=(\arabic*)]
        \item we can run the $K_X+B$ MMP over $U$,
        \item if $K_X+B$ is pseudo-effective, and either $B$ or $K_X+B$ is big over $U$, then any MMP with scaling of a relatively ample divisor terminates with a minimal model, and
        \item if $K_X+B$ is not pseudo-effective over $U$, then any MMP with scaling of a relatively ample divisor terminates with a Mori fiber space.
    \end{enumerate}
    
\end{theorem}

\begin{remark}
    The above theorem allows us to run the MMP only for projective morphisms of analytic varieties. However, the MMP theory for higher-dimensional analytic varieties, even in the case of compact Kähler varieties, is still not well understood. So far, the theory for Kähler 3-folds has been well developed, and \cite{DAS2024109615} established the MMP for Kähler 4-folds in certain specific cases.
\end{remark}

Next, we state two analogs of the negativity lemma in the analytic setting for later use.

\begin{definition}[]
Let $f : X \to Y$ be a proper surjective morphism of normal varieties and $D$ a $\mathbb Q$-divisor. Then
\begin{enumerate}[label=(\arabic*)]
    \item $D$ is $f$-exceptional if $\operatorname{codim}(\operatorname{Supp}(f(D))) \geq 2$.
    \item $D$ is of insufficient fiber type if $\operatorname{codim}(\operatorname{Supp}(f(D))) = 1$ and there exists a prime divisor $\Gamma\not\subset\operatorname{Supp}(D)$ such that $f(\Gamma) \subset \operatorname{Supp}(f (D))$ has codimension one in $Y$.
\end{enumerate}
We say $D$ is degenerate if $D$ is in either of the above cases.
\end{definition}

\begin{lemma}[{\cite[Lemma 2.8]{das2024transcendentalminimalmodelprogram}}]\label{negtrans}
Let $f : X \to Y$ be a proper morphism of normal analytic varieties, where $X$ is a Kähler space, $f_*\mathcal{O}_X = \mathcal{O}_Y$ and $Y$ is 
relatively compact. Let $E = \sum a_iE_i \leq 0$ be a degenerate divisor such that $-E$ is $f$-nef. Then $E = 0$.
\end{lemma}

The following lemma is a relative version of \cite[Lemma 2.9]{lai2011varieties} in the analytic setting.

\begin{lemma}\label{neg}
    Let $f: X\to Y$ be a projective morphism of normal analytic varieties such that $f_*\mathcal{O}_X = \mathcal{O}_Y$ and $D$ be a degenerate effective divisor on $X$. Then there exists a divisorial component of $\operatorname{Supp}(D)$ which is contained in $\mathbf{B}_-(D/Y)$.
\end{lemma}

\begin{pf}
    Let $A$ be a $f$-ample divisor on $X$. We first assume $D$ is of insufficient fiber type. Suppose $Q$ is a divisor in $f(\operatorname{Supp}(D))$ and suppose there exists a divisor on $X$ which dominates $Q$ and is not contained in $\operatorname{Supp}(D)$. Let $F$ be a general fiber for $f$ over $Q$. Working on $F$, we may assume $Q$ is a Cartier divisor. Let $\lambda = \operatorname{min}\{t\geq 0|tf^*Q-D\geq 0\}$. Then $\lambda>0$ and $\lambda f^*Q-D$ is effective. We may take a prime divisors $P\subset\operatorname{Supp}(D)$ and $P'$ over $X$ dominating $Q$ such that $\operatorname{mult}_P(\lambda f^*Q-D) = 0$, $\operatorname{mult}_{P'}(\lambda f^*Q-D) >0$, and $P\cap P'\cap F\neq\emptyset$.
    
    Let $\Bar{F} = F\cap P$ and $\dim \Bar{F} = f$. Then $(A|_{\Bar F})^{f-1}$ represents a general curve $C$ on $\Bar F$. Since $D.C = -(\lambda f^*Q-D).(A|_{\Bar F})<0$, there exists $\epsilon >0$ such that $(D+\epsilon A).C <0$. This implies $C\subset \mathbf{B}_-(D/Y)$, and hence $P\subset\mathbf{B}_-(D/Y)$.
    
    Now, we assume $\operatorname{codim}(f(\operatorname{Supp}(D))\geq 2$. Working locally over a neighborhood of a point $q\in Y$, we may assume $Y$ is Stein. We first suppose $\dim Y = 2$ and $f(\operatorname{Supp}(D)) = \{q\}$ is a point. Let $H$ be a general hyperplane in $Y$ passing through $q$. Let $\lambda = \operatorname{min}\{t\geq 0|tf^*H-D\geq 0\}$. Then $\lambda >0$, $\lambda f^*H-D$ is effective, and there exists a divisorial component $P\subset \operatorname{Supp}(D)$ such that $\operatorname{mult}_{P}(\lambda f^*H-D) = 0.$ By the same argument as above, $D.(A|_{P})^{n-2} = -(\lambda f^*H-D).(A|_{P})^{n-2}<0$, where $n = \dim X$. This implies that $P$ is covered by general curve $C$ with $D.C<0$, therefore, $P\subset \mathbf{B}_-(D/Y).$
    
    Suppose that $\dim Y >2$. Cutting $Y$ by general hyperplane $H_i$ passing through $q$, we produce a surface $Y' = Y\cap H_1\cap \cdots H_n$ containing $q$ and a new map $f':X'=f^{-1}(Y')\to Y'$. By taking normalization, we can further assume $X'$ and $Y'$ are normal. From the surface case we discussed above, there exists a divisorial component in $\operatorname{Supp}(D)\cap X'$ contained in $\mathbf{B}_-(D/Y)$. Since $D$ consists of only finitely many components and $f$, we can find a common divisorial component $P\subset\operatorname{Supp}(D)$ with $P\subset \mathbf{B}_-(D/Y).$
    \hfill$\square$
\end{pf}

Suppose that $f:X\to Y$ is a  projective morphism of compact normal analytic varieties, where $(X,B)$ is klt. Since the divisorial part of $\mathbf{B}_-(K_X+B/Y)$ has only finitely many components, by Theorem~\ref{kahlermmp} and Lemma~\ref{neg}, there exists a model $(X, B)\dashrightarrow{}(X', B')$ over $Y$ such that $\left|K_X'+B'/Y\right|_\mathbb Q$ contains no degenerate divisor.

\subsection{Good Minimal Models}

In this subsection, we state some basic properties of good minimal models, which mostly follow from \cite{BCHM} and \cite{lai2011varieties}. In fact, there is little difference when generalizing these results to the analytic case.

\begin{definition}
    Let $f:X\to U$ be a morphism between normal compact analytic varieties. The pair $(X, B)$ is called a good minimal model over $U$ if $(X,B)$ is a minimal model over $U$ and $K_X+B$ is $f$-semiample.
\end{definition}

\begin{lemma}[cf. {\cite[Lemma 2.1]{lai2011varieties}}]\label{lem1}
    Let $(X, B)$ and $(X', B')$ be compact analytic klt pairs over $U$ and $\mu: X\dashrightarrow{}X'$ a $(K_X+B)$-negative contraction with $\mu_*B = B'$. Then $(X, B)$ has a good minimal model over $U$ if $(X', B')$ does.
\end{lemma}

\begin{pf}
    This follows from \cite[Lemma 3.6.9]{BCHM}.\hfill $\square$
\end{pf}

\begin{lemma}[cf. {\cite[Lemma 2.2]{lai2011varieties}}]\label{lem2}
    Let $(X, B)$ be a compact analytic terminal pair over $U$. Then for any resolution $\mu: (X', B')\to (X, B)$ with $B' = \mu_*^{-1}B$, a good minimal model of $(X', B')$ over $U$ is a good minimal model of $(X, B)$ over $U$.
\end{lemma}

\begin{pf}
    The proof follows from \cite[Lemma 3.6.10]{BCHM}.
    \hfill $\square$
\end{pf}

\begin{theorem}[cf. {\cite[Theorem 2.3, Proposition 2.4]{lai2011varieties}}]\label{uniquegmm}
    Let $f:X\to U$ be a projective morphism between normal compact analytic varieties. Suppose $\phi_i:(X, B)\dashrightarrow{}(X_i, B_i)$ $i = 1,2$ are two minimal models of a klt pair $(X, B)$ over $U$ with $(\phi_i)_*B = B_i$. Then the natural birational map $\psi:(X_1, B_1)\dashrightarrow{}(X_2, B_2)$ over U can be decomposed into a sequence of $(K_{X_1} + B_1)$-flops over $U$. In particular, 
    \begin{enumerate}[label=(\arabic*)]
        \item if $(X, B)$ has a good minimal model over $U$, then any other minimal model over $U$ is also good.
        \item if $(X_1, B_1)$ itself is a good minimal model, then $(X_2, B_2)$ is also a good minimal model.
    \end{enumerate}
\end{theorem}

\begin{pf}
    Since $f$ is a projective morphism, we have the relative cone theorem (\cite[Theorem 2.44]{DAS2024109615}) and we can run relative minimal model programs (\cite[Theorem 1.4]{DAS2024109615}). Therefore, all the arguments in \cite{flop} and \cite[Theorem 2.3]{lai2011varieties} work in our case. (1) follows directly from \cite[Proposition 2.4]{lai2011varieties}. (2) holds because the flops connecting $(X_1, B_1)$ and $(X_2, B_2)$ are crepant.
    \hfill $\square$
\end{pf}

\begin{proposition}[cf. {\cite[Proposition 2.7]{lai2011varieties}}]
\label{GFisGMM}
    Let $f: X\to Y$ be a projective morphism of compact analytic varieties with a general fiber $F$ and $(X,B)$ a klt pair. Suppose $(F, B_F)$ has a good minimal model. Then there exists a $(K_X+B)$-negative birational map $\phi:X\dashrightarrow{} X'$ with a general fiber $F'$ and $B' :=\phi_*^{-1}B$ such that $(F', B'_{F'})$ is a good minimal model.
\end{proposition}

\begin{pf}
     Pick a $f$-ample divisor $H$. By \cite[Theorem 1.4]{DAS2024109615}, we can run a relative MMP with the scaling of $H$. Since the general fiber of $f$ is projective, the proof follows from the rest of argument in \cite[Proposition 2.7]{lai2011varieties}.
    \hfill $\square$
\end{pf}

\subsection{Albanese maps}

We define the Albanese maps for compact Kähler varieties with rational singularities as in  \cite[Definition 2]{das2024existence}.

\begin{definition}
    Let $X$ be a compact Kähler variety with rational singularities $\mu:Y\to X$ a resolution such that $Y$ is a Kähler manifold. Let $a_Y: Y\to \operatorname{Alb}(Y)$ be the Albanese map of $Y$ where $\operatorname{Alb}(Y)$ is a complex torus. Then from the proof of \cite[Lemma 8.1]{Kawamata1985}, $a_Y\circ\mu$ extends uniquely to the morphism $a_X: X\to \operatorname{Alb}(X):=\operatorname{Alb}(Y)$. We call $a_X$ the Albanese map of $X$ and $\operatorname{Alb}(X)$ the Albanese variety of $X$.
\end{definition}

\subsection{Generic Vanishing Theorem}
In this subsection, we recall the following result from the generic vanishing theory in the Kähler setting.

\begin{definition}
    Let $X$ be a complex manifold and $\mathcal F$ a coherent sheaf on $X$. Then the $k$-th jumping locus of the $i$-th cohomology is defined by
    $$V^i_k(\mathcal F) = \{\rho\in\operatorname{Pic}^0(X)|h^i(X, \mathcal F\otimes \rho)\geq k\}.$$
\end{definition}

\begin{theorem}[{\cite[Corollary 1.4]{wang2016torsion}}]\label{GV}
    Let $X$ be a compact Kähler manifold. Then $V^i_k(K_X)$ is a finite union of torsion translate of subtori in $\operatorname{Pic}^\circ (X)$.
\end{theorem}

\begin{remark}
    In \cite[Theorem D]{AFST_2021_6_30_4_813_0}, the above theorem is settled in a more general setting.
\end{remark}

\section{Proof of the Main Theorems}
First, we will show that \( \kappa(X, K_X + B) \geq 0 \) under the assumptions of Theorem~\ref{1}. The following theorem demonstrates that \cite[Theorem A (II)]{AFST_2021_6_30_4_813_0} applies to any map to a complex torus with connected fibers, even if the map is not surjective.

\begin{theorem}\label{nonvanish}
    Let $a: X\to A$ be a morphism with connected fibers from a compact Kähler variety $X$ to a complex torus $A$. Let $F$ be the general fiber of $a:X\to a(X)=Z$, and $(X, B)$ a klt pair. Then
    $$\kappa(X, K_X+B)\geq\kappa(F, K_F+B_F) + \kappa(Z).$$ 
\end{theorem}

\begin{pf}
    Passing to a terminalization by Theorem~\ref{kahlermmp} and \cite[Corollary 1.4.3]{BCHM}, we may assume $(X, B)$ has terminal singularities. Let $Z^\nu$ be the normalization of $Z$, which induces a surjective morphism $a^\nu: X\to Z^\nu$. 
    
    Apply \cite[Proposition 4.1]{AFST_2021_6_30_4_813_0} to $Z^\nu\to A$, there exists a proper surjective morphism $\phi_{a^\nu}:Z^\nu\to W$ such that $W$ is a normal variety of general type, and the general fiber is a complex torus $\tilde S$. Let $X_w$ be the general fiber of $X\to W$. By \cite[Theorem 3.2]{AFST_2021_6_30_4_813_0},
    \begin{align*}
        \kappa(X, K_X+B) &\ge \dim W + \kappa(X_w, K_{X_w} + B_{X_w})\\
        &= \kappa(Z^\nu) + \kappa(X_w, K_{X_w} + B_{X_w}).
    \end{align*}

    On the other hand, $a^\nu$ maps $X_w$ onto $\tilde S$. Then, by \cite[Theorem A (II)]{AFST_2021_6_30_4_813_0},
    $$
    \kappa(X_w, K_{X_w} + B_{X_w}) \ge \kappa(F, K_F+B_F) + \kappa(\tilde S) \ge \kappa(F, K_F+B_F),
    $$
    which completes the proof.
    \hfill $\square$
\end{pf}

Next, the following lemma will help prove Theorem~\ref{1} in the case where $\kappa(X, B) = 0$. In \cite[Lemma 4.1]{lai2011varieties}, Lai has proved a similar statement without considering log pairs. To tackle the case of log pairs, we first need to take a suitable cyclic covering, apply a similar argument to the covering space, and then transfer the information back to $X$.

\begin{lemma}\label{insuf}
    Let $X$ be a compact Kähler manifold, $(X, B)$ a klt pair such that $\kappa(X, K_X+B) = 0$, and $a_X: X\to A$ be the Albanese map. Suppose that $\Gamma\in |m(K_X+B)|$ is the unique effective divisor, then $\operatorname{Supp}(\Gamma)$ contains all $a_X$-exceptional divisors.
\end{lemma}

\begin{pf}
    By taking a resolution of $(X, B+\Gamma)$, we may assume that $B+\Gamma$ has an SNC support. Also, we may assume that  $mB$ is an integral divisor and $B$ and $\Gamma$ have no common components, otherwise, we simply subtract the common divisor from $\Gamma$ and $B$.
    
    Let $D = \Gamma + m(\lceil B\rceil - B)\in |m(K_X+\lceil B\rceil)|$. We take a normalized cyclic covering along $D$ (cf. \cite[Lemma 1.1]{AFST_2021_6_30_4_813_0} and \cite[Proposition 2.1]{Kawamata1985}), say $\pi: V\to X$, such that $V$ is smooth and for all $p\geq 0$,
    $$\pi_*\Omega_V^p = \bigoplus_{i = 0}^{m-1}\Omega_X^p(\log D_i)\otimes\mathcal{O}_X\left({-i(K_X+\lceil B\rceil) + \left\lfloor\frac{i}{m}D\right\rfloor}\right),$$
    where $D_i  = \left\{\frac{i}{m}D\right\}_{red}$. In particular,
    $$\pi_*\omega_V = \bigoplus_{i = 0}^{m-1}\omega_X( D_i)\otimes\mathcal{O}_X\left({-i(K_X+\lceil B\rceil) + \left\lfloor\frac{i}{m}D\right\rfloor}\right).$$
    The $i = 1$ summand of $\pi_* \omega_V$ is 
    \begin{align*}
            &K_X+\left\{\frac{D}{m}\right\}_{red}-\left(K_X+\lceil B\rceil\right) +\left\lfloor\frac{D}{m}\right\rfloor\\
            =&K_X+\left\{\frac{\Gamma}{m}\right\}_{red} + \lceil B\rceil - \left(K_X+\lceil B\rceil\right) + \left\lfloor\frac{\Gamma}{m}\right\rfloor = \left\lceil\frac{\Gamma}{m}\right\rceil.
    \end{align*}
    We denote $L = \mathcal O_X\left(\left\lceil\frac{\Gamma}{m}\right\rceil\right)$. Then $H^0(X, L)\neq 0$ and $\kappa(L) = 0$. By \cite[Lemma 5.6]{AFST_2021_6_30_4_813_0}, each irreducible component of $V^0_1(L)$ is also an irreducible component of $V^0_k(\pi_*\omega_V)$ for some $k\geq 1$. Then, $V^0_1(L)$ is a finite union of torsion translates of subtori in $\operatorname{Pic}^\circ(X)$ by Theorem~\ref{GV}. Following the proof of \cite[Proposition 2.1]{ein1997singularities}, one can show that $0$ is an isolated point in $V^0_1(L)$.

    Now, let $v\in H^1(X, \mathcal O_X)$ and $\eta = \bar v\in H^0(X, \Omega_X^1)$ is its conjugate. We can define the following derivative complex
    \begin{equation}
        H^0(V, \Omega_V^{n-1})\xrightarrow{\wedge\pi^*\eta} H^0(V, \omega_V)\xrightarrow{\wedge\pi^*\eta} 0.
    \end{equation}
    Applying the above decomposition, we claim that the complex on the $i = 1$ summand,
    \begin{equation}
        H^0\left(X, \Omega_X^{n-1}(\log D_1)\otimes\mathcal{O}_X\left({-(K_X+\lceil B\rceil) + \left\lfloor\frac{D}{m}\right\rfloor}\right)\right)\xrightarrow{\wedge \eta}H^0(X, L)\xrightarrow{\wedge \eta}0,
    \end{equation}
    is exact for every $v\neq 0$. The argument is similar to the proof of \cite[Theorem 1.2.3]{ein1997singularities}. Let $\Delta_{v}(0)\subset\operatorname{Pic}^0(X)$ be a neighborhood of $0$ in the straight line in $\operatorname{Pic}^0(X)$ through $0$ determined by $v$. By \cite[Corollary 3.3]{green1991higher}, when $\Delta_v(0)$ is sufficiently small, the homology of (1) is isomorphic to 
    $$H^0(V, \omega_V\otimes \pi^*\rho) = H^0(X, \pi_*\omega_V\otimes \rho),$$ 
    for $\rho\in\Delta_v(0)$. This implies the homology of (2) is isomorphic to $H^0(X,L\otimes \rho) = 0$, since $0$ is an isolated point in $V^0_1(L)$, which proves the claim.

    Now, given an $a_X$-exceptional prime divisor $E\subset X$. Following the first paragraph of the proof of \cite[Lemma 4.1]{lai2011varieties}, there exists a nonzero $1$-form $\eta_e$ which is given by pulling back of a flat $1$-form on $A$ and vanishes on at the general point $e\in E$. Then by the surjectivity of (2), $E$ is contained in $\operatorname{Supp}(L)$ = $\operatorname{Supp}(\Gamma)$, which completes the proof.
\hfill $\square$

\end{pf}

\begin{remark}
    In fact, Lemma~\ref{insuf} holds for any morphism from $X$ to a complex torus. All the arguments in the proof remain valid as long as we work with $\operatorname{Pic}^0(A)$ instead of $\operatorname{Pic}^0(X)$. 
\end{remark}

Before proving Theorem~\ref{1}, we establish the following theorem, which can be viewed as a replacement for the cone theorem in the projective case. For generalized pairs over Kähler varieties, we refer the reader to \cite{hacon2024canonicalbundleformulaadjunction}.

\begin{theorem}\label{main1}
    Let $a:X\to A$ be a projective morphism with connected fiber from a compact Kähler variety $X$ to a complex torus $A$, and  $(X, B)$ a klt pair. Let $F$ be the general fiber of $a:X\to a(X)$.
    Suppose $(X, B)$ is a relatively minimal model over $A$ and $\kappa(F, K_F+B_F)\geq0$. Then $(X, B)$ is a minimal model.
\end{theorem}

\begin{pf}
    Let $H$ be an $a$-ample divisor and $\omega_A$ be a Kähler form on $A$ such that $w := H + a^*\omega_A$ is a Kähler form on $X$. We assume $(X, B)$ is not a minimal model. Then we can find $t>0$ such that $\alpha:=K_X+B+tw$ is nef but not Kähler. By Theorem~\ref{nonvanish}, $K_X+B$ is effective, therefore, $\alpha$ is also big.
    
    By \cite[Theorem 1.1]{Collins_2015}, we can take $\eta$ in the class $\alpha$ as a Kähler current with analytic singularities such that $\operatorname{Null}(\alpha) = E_{+}(\eta)$, where $\operatorname{Null}(\alpha)\neq \emptyset$. We let $Z$ be a maximal dimensional irreducible component of $\operatorname{Null}(\alpha)$.

    We first follow the argument in \cite[Proposition 3.1]{hacon2024canonicalbundleformulaadjunction}. Let $c$ be the log canonical threshold of $(X, B)$ with respect to $\eta$ near the general point of $Z$. We consider a resolution $\nu: X'\to X$ such that $K_{X'}+B' = \nu^*(K_X+B)$, $\nu^*\eta = \eta'+F$ where $\eta'$ is semi-positive and $F$ is an effective $\mathbb R$-divisor, and $F+B'$ has a simple normal crossing support. We define the $b$-$(1,1)$-form $\boldsymbol{\eta} = \overline{\eta'}$. Then near the general point of $Z$, $(X, B+c\nu_*F+c\boldsymbol{\eta})$ is a generalized log canonical pair and $Z$ is its minimal log canonical center. By \cite[Remark 2.12]{hacon2024canonicalbundleformulaadjunction}, we have
    $$(K_X+B+c\nu_*F+c\boldsymbol{\eta})|_{Z^\nu}= K_{Z^\nu}+ B_{Z^\nu}+ \boldsymbol{\beta}_{Z^\nu},$$
    where $Z^\nu$ is the normalization of $Z$, and $(Z^\nu, B_{Z^\nu}+\boldsymbol{\beta}_{Z^\nu})$ is a generalized pair over $Z^\nu$.

    Now, consider the map $Z^\nu\to A$ induced by $a$. Let $V$ be the Stein factorization of $Z^\nu\to A$, $V'$ its resolution, and $Z'$ a resolution of $(Z, B_{Z^\nu}+\boldsymbol{\beta}_{Z^\nu})$ such that
    $$\begin{tikzcd}
    Z' \arrow{d}{\phi'}\arrow{r}{\mu}&Z^\nu \arrow{d}{\phi}{}\arrow{r}{a}  & A \\
    V'\arrow{r}&V\arrow{ru}{}
    \end{tikzcd}$$
    commutes. Denote the general fiber of $\phi'$ and $\phi$ by $G'$ and $G$ respectively. We write $K_{Z'}+B_{Z'}+\boldsymbol{}\beta_{Z'} = \mu^*(K_{Z^\nu} + B_{Z^\nu}+\beta_{Z^\nu}) + E,$
    where $B_{Z'} = \mu^{-1}_*B_{Z^\nu}$, $\beta_{Z'} = \mu^*\beta_{Z^\nu}$, and $E$ is an exceptional divisor. Since $H$ is $a$-ample, we can take a sufficiently large $l$ such that $h^0(G', E|_{G'}+lH|_{G'})>0$. We now consider the klt generalized pair $(Z', B_{Z'} + \beta_{Z'}+\delta\eta_{Z'})$, where $\boldsymbol{\eta}_{Z'} = \boldsymbol{\eta}|_{Z'}$. We have
    \begin{align*}
        &K_{Z'} + B_{Z'} + \beta_{Z'} +\delta\eta_{Z'}\\ 
        =&((c+\delta+1)(K_X+B)+ctH)|_{Z'} + E + \delta tH|_{Z'}+((c+\delta)t)\mu^*a^*\omega_A.
    \end{align*} 
    Since $(c+1)(K_X+B)+ctH$ is $a$-ample, for sufficiently large $m$,
    $$h^0(G', m((c+1)(K_X+B)+ctH))>0,$$
    and hence,
    $$h^0(G', m(((c+\delta+1)(K_X+B)+ctH))|_{Z'} + E|_{G'}+\delta tH|_{G'}))>0$$
    for sufficiently large $\delta$. By applying \cite[Theorem 2.2]{hacon2024canonicalbundleformulaadjunction} to $\phi'$, $K_{Z'/V'} + B_{Z'} + \beta_{Z'}+\delta \eta_{Z'}$ is pseudoeffective.  On the other hand, $\kappa(V')\geq 0$ by \cite[Proposition 4.1]{AFST_2021_6_30_4_813_0}. Thus, we conclude that $K_{Z'} + B_{Z'} + \beta_{Z'}+\delta \eta_{Z'}$ is pseudoeffective.

    Let $k<\dim Z$ be the numerical dimension of $\alpha_{Z^{\nu}}$, that is, $\alpha_{Z^\nu}^k\neq 0$ and $\alpha_{Z^\nu}^{k+1} = 0$. Then 
    \begin{align*}
        &(K_{Z'}+B_{Z'} + \beta_{Z'}+\delta \eta_{Z'}).\alpha^{k}_{Z'}\cdot\omega_{Z'}^{\dim Z-k-1}\\
        =&(K_{Z^\nu}+B_{Z^\nu}+\beta_{Z^\nu})\cdot\alpha_{Z^\nu}^k\cdot\omega_{Z^\nu}^{\dim Z-k-1}\\
        =&(K_X+B+c\alpha)|_{Z^\nu}\cdot\alpha_{Z^\nu}^k\cdot\omega_{Z^\nu}^{\dim Z-k-1}\\
        =&((c+1)\alpha-t\omega)|_{Z^\nu}\cdot\alpha_{Z^\nu}^k\cdot\omega_{Z^\nu}^{\dim Z-k-1}\\
        =&-t \alpha_{Z^\nu}^k\cdot\omega_{Z^\nu}^{\dim Z-k}<0.
    \end{align*}
    This implies $K_{Z'} + B_{Z'} + \beta_{Z'}+\delta \eta_{Z'}$ cannot be pseudoeffective, a contradiction, which concludes the proof.\hfill$\square$
\end{pf}

\begin{remark}
In \cite[Theorem~3.1]{hacon2024canonicalbundleformulaadjunction}, Hacon and P\v{a}un establish an analytic cone theorem for $n$-dimensional Kähler klt generalized pairs, assuming the BDPP conjecture \cite[Conjecture~0.1]{boucksom2004pseudo} holds in dimension
 $\le n-1$. As a consequence, one obtains an $\alpha$-trivial rational curve on $X$ which cannot be contracted by the morphism $a$. This leads to a contradiction as in Fujino's argument mentioned at the beginning of this section.

Let $X$ be a Kähler manifold. The BDPP conjecture asserts that $K_X$ is pseudoeffective if and only if $X$ is not uniruled, that is, not covered by rational curves. At the time of writing this thesis, this conjecture has been proved in \cite[Theorem~1.1]{ou2025characterization}.
\end{remark}

\begin{pf2}
    By applying the $\mathbb{Q}$-factorial terminalization (via Theorem~\ref{kahlermmp} and \cite[Corollary 1.4.3]{BCHM}), we can assume that $(X, B)$ is both $\mathbb{Q}$-factorial and terminal. Additionally, using Lemma~\ref{lem2}, we can further assume that $(X, B)$ is log smooth. According to Proposition~\ref{GFisGMM} and Lemma~\ref{lem1}, we may also assume that $(F, B_F)$ is a good minimal model. Finally, by Theorem~\ref{nonvanish}, we conclude that $\kappa(X, K_X + B) \geq 0$.

    We first address the case where $\kappa(X, K_X + B) = 0$, which is quite similar to the projective case (\cite[Theorem 4.2]{lai2011varieties}). When $\kappa(X, K_X+B) = 0$, $a_X:X\to A$ is a fiber space (\cite[Theorem C]{AFST_2021_6_30_4_813_0}). By running a relative MMP with scaling, we may further assume that $\operatorname{\bf{B}}_-
    (K_X+B/A)$ contains no divisorial component. 
    Using \cite[Theorem A]{AFST_2021_6_30_4_813_0}, we obtain $\kappa(F, K_F+B_F)\leq\kappa(X, K_X+B) = 0$. Since $K_F+B_F$ is semiample, we conclude that $\kappa(F, K_F+B_F) = 0$ and $K_F+B_F\sim_{\mathbb{Q}}\mathcal{O}_F$.

    Let $\Gamma\geq 0$ with $K_X+B\sim_{\mathbb Q}\Gamma$. We aim to prove $\Gamma = 0$. From the argument in \cite[Theorem 4.2]{lai2011varieties} as well as Lemma~\ref{insuf}, it follows that $\Gamma$ is degenerate. If $\Gamma\neq0$, then by Lemma~\ref{neg}, there exists a divisorial component of $\operatorname{Supp}(\Gamma)$ which is contained in $\mathbf{B}_-(K_X+B/A)$, which leads to a contradiction. Therefore, we conclude that $\Gamma = 0$ and $K_X+B$ is semiample.

    Now, we assume $\kappa(X, K_X+B) > 0$. By Theorem~\ref{fg}, $\displaystyle R(X, K_X+B)$ is finitely generated. We fix $d\in\mathbb N$ such that $R^{[d]}(X, K_X+B)$ is generated in degree $1$. Let $\mu: X'\to X$ be a resolution of $|d(K_X+B)|$. Then, we have
    \begin{center}
        \begin{tikzcd}
            X'\arrow[d, "\mu"] \arrow[rd, "g"] \\
            X \arrow[r, dashed]& Y = \operatorname{Proj} R^{[d]}(X,K_X+B),
        \end{tikzcd} 
    \end{center}
    where $X\dashrightarrow{} Y$ is the Iitaka fibration.

    Let $B'=\mu_*^{-1}B$ and $X'_y$ be the general fiber of $g$. Then $\kappa(X'_y, K_{X'_y}+B'_{X'_y}) = 0$. We will first show that $(X'_y, B'_{X'_y})$ on the general fiber of the Albanese map of $X'_y$ has a good minimal model, which, by the result from the previous case, implies that $(X'_y, B'_{X'_y})$ itself also has a good minimal model. By \cite[Proposition 8.]{das2024existence}, $Y$ has rational singularities, so we can consider the following diagram.
    \begin{center}
        \begin{tikzcd}
            X'_y\arrow[r]\arrow[d, "a_{X'_y}"]&X'\arrow[d, "a_{X'}"] \arrow[r, "g"]  &Y\arrow[d, "a_{Y}"]\\
            A(X'_y)\arrow[r]&A = A(X)\arrow[r] &A(Y),
        \end{tikzcd} 
    \end{center}
    where $a_{X'_y}$, $a_{X'}$, and $a_Y$ are Albanese maps. Note that $a_{X'}$ factors through $f$, so $(X', B')$ is a good minimal model on the general fiber of $a_{X'}$. 
    
    By \cite[Lemma 2.6]{doi:10.1081/AGB-120027861}, $a_{X'}(X'_y)$ is a translation of a fixed subtorus $K\subset A(X)$ for general $y\in Y$, and $0\to K\to A(X)\to A(Y)\to 0$ is exact. Then we have $(X'_y, B'_{X'_y})$ is a good minimal model on the general fiber of $a_{X'_y}$, which proves $(X'_y, B'_{X'_y})$ also has a good minimal model.

    Let $n = \dim X'$ and let $M' = g^*\mathcal O_Y(1)$ be the base point free divisor defining $g$. We can take an effective $\mathbb Q$-divisor $H'$ such that $H'\sim_{\mathbb Q}(2n+1)M'$ and $(X',B'+H')$ is still klt. In particular, $(X', B'+H')$ is a good minimal model on the general fiber of $a_{X'}$. Follow the proof of \cite[Theorem 2.12]{hacon2013existence} and Theorem~\ref{kahlermmp}, by running a $(K_{X'} + B'+H')$-MMP with a scaling over $A$, we obtain $\phi:X'\dashrightarrow{} X''$ over $A$, a relatively good minimal model of $(X', B'+H')$. By Theorem~\ref{main1}, $K_{X''}+B''+H''$ is nef.
    
    We will show $(X'', B''+H'')$ is also a model over $Y$. Let $R\in \operatorname{\overline{NE}}(X'/A)$ be a $(K_{X}+B'+H)$-negative extremal ray, then it is also a $(K_{X'} + B')$-negative extremal ray. By the relative cone theorem \cite[Theorem 2.45]{DAS2024109615}, $R$ is spanned by a rational curve $C$ such that $0 > (K_{X'} + B').C \geq-2n$. On the other hand, if $g(C)$ is not a point, $C.M'\geq1$. Then 
    $$(K_{X'} + B' + H').C = (K_{X'} + B' + (2n+1)M').C>0,$$
    which leads to a contradiction. Thus, $C$ must be contracted by $g$ and $X''$ is a model over $Y$. We have the following diagram.
    \begin{center}
        \begin{tikzcd}
            X'\arrow[d, "g"]\arrow[rd, ]\arrow[r, dashed, "\phi"]&X''\arrow[d,]\arrow[ld,"h"]\\
            Y&A.
        \end{tikzcd} 
    \end{center}
    Let $M''=h^*\mathcal O_Y(1)$ be the divisor defining $h$. Then for every fiber $X''_y$ over $Y$,
    $$(K_{X''} + B''_{X''})|_{X''_y} = (K_{X''} + B''_{X''}+(2n+1)M'')|_{X''_y} = (K_{X''} + B''_{X''} + H'')|_{X''_y}$$
    is nef. Thus $K_{X''}+B_{X''}$ is nef over $Y$. 
    
    We will show that for the general fiber $X''_y$ of $X''\to Y$, $(X''_y, B''_{X''_y})$ is actually a good minimal mode. Indeed, the Albanese map of $X_y$ is also projective. Since $(X'_y, B_{X'_y})$ has a good minimal model and $(X''_y, B''_{X''_y})$ is a minimal model, $(X''_y, B''_{X''_y})$ is actually good by (2) of Theorem~\ref{uniquegmm}. This implies $K_{X''_y} + B''_{X''_y}\sim_\mathbb Q 0$.

    On the other hand, we write $d\mu^*(K_X+B) = M' + F'$, where $M'$ is the free part and $F'$ is the fixed part. And suppose $K_{X'} +B' = \mu^*
    (K_X+B) + E'$, where $E$ is an effective divisor since $(X, B)$ is terminal. Then we have $d(K_{X'} + B') = M' + dF' + E'$. Pushing forward on $X''$, we have $d(K_{X''} + B'') = M'' + dF''+E''$. By restricting it to the general fiber,
    $$(dF''+E'')|_{X''_y} = M'' + dF''+E''|_{X''_y} =d(K_{X''} + B'')_{X''_y}\sim _\mathbb Q 0.$$
    This implies $dF''+E''$ has no horizontal part. In addition, from the argument in \cite[Theorem 4.4]{lai2011varieties}, $dF''+E''$ is degenerate. Then we conclude $dF''+E'' = 0$ by Lemma~\ref{negtrans}. Thus, $(X'', B'')$ is a good minimal model, which implies $(X, B)$ has a minimal model by Lemma~\ref{lem2}.
    \hfill $\square$
\end{pf2}

\begin{remark}
By the universal property of the Albanese map, in Theorem~\ref{1} and Lemma~\ref{insuf} it is not necessary to assume that $a_X$ is the Albanese map. Both results remain valid for any projective morphism to a complex torus.
\end{remark}

\section*{Acknowledgments}
The author would like to thank her advisor, Christopher Hacon, for suggesting this problem and for his helpful guidance and discussions throughout this work.

The author was partially supported by NSF research grants no. DMS-1952522 and DMS-1801851, as well as by a grant from the Simons Foundation (SFI-MPS-MOV-00006719-07).

\bibliographystyle{alpha}  
\bibliography{main1}

\end{document}